\documentclass[12pt,thmsa]{article}

\usepackage[ansinew]{inputenc}
\usepackage{pst-all}
\usepackage{amsmath}
\usepackage{amsfonts}
\usepackage{amssymb}

\newtheorem{lemma}{Lemma}
\newtheorem{proposition}[lemma]{Proposition}
\newtheorem{theorem}[lemma]{Theorem}
\newtheorem{example}[lemma]{Example}
\newtheorem{remark}[lemma]{Remark}

\newcommand{\fin}{\hfill $\Box$}

\title {The equal tangents property}
\author{J. Jer\'onimo-Castro, G. Ruiz-Hern\'andez and S. Tabachnikov}

\begin{document}
\maketitle

\begin{abstract}
Let $M$ be a $C^2$-smooth strictly convex closed surface in
$\mathbb{R}^3$ and denote by $H$ the set of points $x$ in the
exterior of $M$ such that all the tangent segments from $x$ to $M$
have equal lengths. In this note we prove that if $H$ is either a
closed surface containing $M$ or a plane, then $M$ is an Euclidean
sphere. Moreover, we shall see that the situation in the Euclidean
plane is very different.
\end{abstract}

\section{In the Euclidean plane} \label{plane}
Let $K$ be a strictly convex body in the plane. The following fact is well
known: \emph{if the two tangent segments to $K$ from
every point $x\not\in K$  have equal lengths then $K$ is an Euclidean disc}
(see, for instance, \cite{RaTop} ). This statement is easily proved by
elementary geometry. The  result was extended to the case of
Minkowski planes by S. Wu \cite{Wu}, and Z. L\'angi \cite{Langi} also gave a
characterization of the ellipsoid
among centrally symmetric convex bodies in terms  of tangent
segments of equal Minkowski length.

In the Euclidean plane, one may obtain the same conclusion with considerably weaker assumptions.
 Namely, one has the following characterization of a circle in terms  of equal tangent segments.

\begin{lemma} \label{line}
Let $\gamma$ be a strictly convex closed curve
in the plane, and let $\ell$ be a tangent line through a point
$p\in \gamma$. Suppose that the two tangent segments to $\gamma$
from every point $x\in \ell$ have equal lengths. Then $\gamma $ is a circle
(see, e.g., \cite{JR} or Section \ref{further}).
\end{lemma}

Thus it is natural to ask whether the same conclusion remains
true if the locus of points from which the tangent segments to $\gamma$ have equal lengths is a
line  $\ell$ that is not tangent to $\gamma$. We consider the two cases separately: first, when
$\ell$ intersects $\gamma$, and second, when $\ell$ is disjoint from $\gamma$.

\begin{example}
{\rm  Let us construct a non-circular curve with the desired equitangent property.}
\end{example}

Consider two circles and their radical axis
$\ell$ (that is, the set of points having equal power with respect
to both of them). Then the tangent segments to both circles from the points of $\ell$ are equal.
Let $x,y\in \ell$ be two points in the
exterior of the convex hull of these circles. Draw the tangents
$xa,xb,yc,$ and $yd$ as shown in Figure 1, and also draw the two
arcs of the circles tangent to $xa,xb$ at $a$ and $b,$ and to
$yc,yd$ at $c$ and $d,$ respectively. Then the
union of the arcs $\widehat{ab},\widehat{bd},\widehat{dc},$ and
$\widehat{ca}$ is a $C^1$-smooth and strictly convex curve $\gamma$ with
the property that for every point $p\in
\ell\setminus\text{conv}(\gamma)$, the two tangent segments to $\gamma$
from $p$ have equal lengths.

\centerline{ \psset{unit=.8cm}
\begin{pspicture}(-1,0)(7.4,7.2)
{\psarc[linecolor=blue,linewidth=.035](2.41,5.65){.96}{53.69}{146.6}
\psarc[linestyle=dashed](2.41,5.65){.96}{146.6}{53.69}
\psarc[linecolor=blue,linewidth=.035](2.41,2.19){1.72}{189.72}{-35.64}
\psarc[linestyle=dashed](2.41,2.19){1.72}{-35.64}{189.72}
\psarc[linecolor=blue,linewidth=.035](6.58,2.91){5.95}{146.6}{189.72}
\psarc[linecolor=blue,linewidth=.035](.75,3.39){3.77}{-35.64}{53.69}
\psline[linecolor=gray,linestyle=dashed](1.61,6.18)(.32,4.22)(.72,1.9)
\psline[linecolor=gray,linestyle=dashed](2.98,6.43)(5.98,4.22)(3.81,1.19)
\psline[linecolor=red](-.74,4.22)(7.14,4.22)
\psdots[dotsize=3pt](1.61,6.18)(.72,1.9)(2.98,6.43)(3.81,1.19)(.32,4.22)(5.98,4.22)
\rput(5.2,4.5){\small $\ell$}\rput(.45,1.8){\small
$b$}\rput(3.1,6.55){\small $c$}\rput(4,1.1){\small $d$}
\rput(1.4,6.3){\small $a$}\rput(.1,4.5){\small
$x$}\rput(6.1,4.5){\small $y$}}\rput(1.5,.4){\small $\gamma$}
\end{pspicture}}
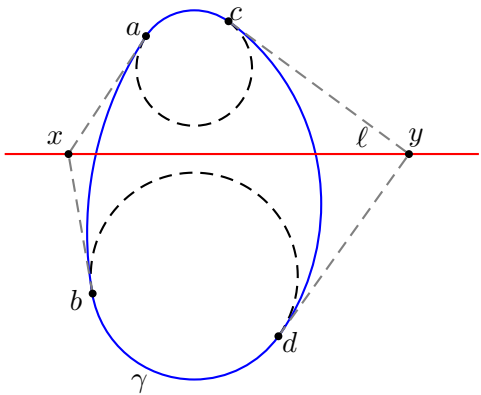
\centerline{\bf\small{Figure 1}}

\textbf{A characterization in terms of hyperbolic geometry.} In the case when $\ell$ does
not intersect $\gamma$, we have a complete characterization of such curves.

\begin{lemma} \label{hyper}
Assume that $\ell$ is the horizontal axis and that $\gamma$ lies in the upper half plane.
Then the tangent segments to $\gamma$ from every point of $\ell$ are equal if and only if
$\gamma$ is a curve of constant width in the hyperbolic metric, considered in the upper half plane model.
\end{lemma}

\emph{Proof.} The two tangent segments from a point  $x\in\ell$ have equal lengths if and
only if the circle centered at $x$ is orthogonal to $\gamma$ at both intersection points
(more precisely, is orthogonal to support lines to $\gamma$ at these points). See Figure 2.

\centerline{ \psset{unit=1.2cm}
\begin{pspicture}(0,.4)(4,3.6)
\psline[linecolor=gray,linestyle=dashed](1.45,2.89)(1.04,.79)(2.91,1.28)
\psline[linecolor=gray,linestyle=dashed](1.78,1.24)(4.14,.79)(2.67,2.69)
\pspolygon[linecolor=lightgray,fillstyle=solid,fillcolor=lightgray](1.56,3.27)(1.78,1.24)(2.91,1.28)
\psarc[linewidth=.035,fillstyle=solid,fillcolor=lightgray](3.60,2.46){2.19}{158}{214}
\psarc[linewidth=.035,fillstyle=solid,fillcolor=lightgray](1.48,1.77){1.51}{-18}{87}
\psarc[linewidth=.035,fillstyle=solid,fillcolor=lightgray](2.26,3.73){2.53}{260}{285}
\psdots[dotsize=3pt](1.56,3.27)(1.78,1.24)(2.91,1.28)(1.04,.79)(4.14,.79)(2.67,2.69)(1.45,2.89)
\psline[linewidth=.03](0,.79)(5.5,.79) \rput(1.5,3.5){\small
$a$}\rput(1.6,1.15){\small $b$}\rput(3.1,1.2){\small $c$}
\rput(1,.5){\small $x$}\rput(4.2,.5){\small
$y$}\rput(5.3,1){\small $\ell$}
\end{pspicture}}
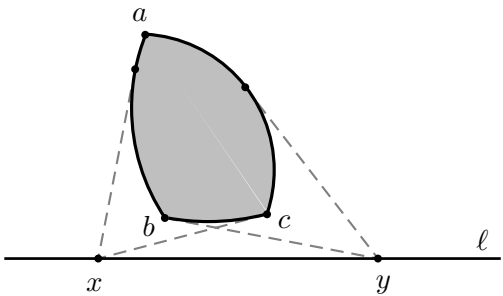
\centerline{\small Figure 2.}

The circles centered at points of $\ell$ are the geodesics of the hyperbolic plane,
and the upper half plane model is conformal. Thus a geodesic segment can make a full
 circuit inside $\gamma$, remaining orthogonal to it at both end points. This
property  characterizes convex bodies of constant  width, see \cite{Ara}.
\fin

Next we construct pairs of curves in the plane, $\Gamma$ and
$\gamma$, such that $\Gamma$ encloses $\gamma$ and, for every
point $x\in \Gamma$, the two tangent segments from $x$ to $\gamma$
have equal lengths. Compare with \cite{Tab3} where a pair of
curves $\Gamma$ and $\gamma$ is constructed such that, for every
point $x\in \Gamma$, the tangent segments from $x$ to $\Gamma$
have unequal lengths.

\begin{example}
{\rm Consider a regular convex $n$-gon, with odd $n\geq 5$, and
make the classical construction of a body of constant width.
Concretely, let $V_1V_2V_3V_4V_5$ be a regular pentagon and let
$\lambda$ be the length of its diagonals. Fix $\varepsilon\geq 0$.
Draw the lines $V_1V_3,\, V_1V_4,\, V_5V_2,\, V_5V_3,$ and
$V_2V_4$. Now, draw the arcs of the circles centered at $V_1$ and
the radii $\varepsilon$ and $\lambda+\varepsilon$  from line
$V_3V_1$ to line $V_4V_1$, see Figure 3. Do the same at the
remaining vertices. We obtain a $C^1$-smooth convex curve $\gamma$
of constant width $\lambda+2\varepsilon$.

\centerline{ \psset{unit=.9cm}
\begin{pspicture}(0,-.5)(7,7)
\psline(3.31,4.58)(2.06,3.66)(2.54,2.19)(4.09,2.19)(4.57,3.66)(3.31,4.58)
\def\chi{\psarc[linewidth=.03,linecolor=blue](0,0){.56}{72}{108}
\psline[linecolor=gray,linestyle=dashed](.17,.52)(-.94,-2.92)
\psline[linecolor=gray,linestyle=dashed](.94,-2.92)(-.17,.52)
\psarc[linewidth=.03,linecolor=blue](.78,-2.39){3.08}{108}{144}
\psline(0,1.75)(-1.8,1.16)(-2.92,-.37)
\psdots[dotsize=3pt](0,1.75)(-1.8,1.16)} \rput{0}(3.31,4.58){\chi}
\rput{72}(2.06,3.66){\chi}
\rput{144}(2.54,2.19){\chi}\rput{216}(4.09,2.19){\chi}\rput{288}(4.57,3.66){\chi}
\psline[linestyle=dashed,linecolor=red](3.32,5.16)(1.09,5.17)(1.55,2.99)
\psdots[dotsize=3pt](3.32,5.16)(1.09,5.17)(1.55,2.99)
\rput(3.6,4.65){\small $V_1$}\rput(1.85,3.4){\small
$V_2$}\rput(2.2,2.35){\small $V_3$}\rput(4.4,2.35){\small
$V_4$}\rput(4.6,4){\small $V_5$} \rput(1.4,6){\small
$P_1$}\rput(.05,4.3){\small $P_2$}\rput(.2,2.1){\small
$P_3$}\rput(1.3,.5){\small $P_4$}\rput(3.5,-.1){\small
$P_5$}\rput(5.3,.5){\small $P_6$}\rput(6.5,2.1){\small
$P_7$}\rput(6.4,4.4){\small $P_8$}\rput(5.3,6){\small
$P_9$}\rput(3.4,6.6){\small $P_{10}$} \rput(.9,5.25){\small
$X$}\rput(5.15,2.4){\small $\gamma$}\rput(5.9,5){\small $\Gamma$}
\end{pspicture}}
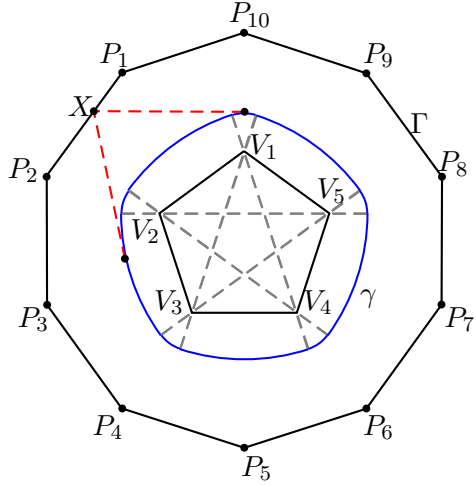
\centerline{\small Figure 3.}

Let $\Gamma =P_1P_2\ldots P_{10}$ be the regular decagon
constructed in the following way: the segment $P_1P_2$ is
contained in the radical axis of the circles with centers $V_1$
and $V_5$ and the radii $\varepsilon$ and $\lambda+\varepsilon,$
respectively; this axis is orthogonal to the side $V_1 V_5$.
Likewise, the segment $P_2P_3$ is contained in the radical axis of
the circles with centers $V_4$ and $V_3$ and the radii
$\varepsilon$ and $\lambda+\varepsilon,$ respectively, etc. Then
the tangent segments to $\gamma$ from every point of $\Gamma$ are
equal.}
\end{example}

\begin{remark}
{\rm It is interesting to investigate what happens if one imposes additional assumptions
on the curve $\Gamma$. For example, is it true that if $\Gamma$ is a circle then $\gamma$ also
must be a circle? We do not know the answer to this question.

We remark that the existence of planar bodies (different from the circle) floating in equilibrium in all positions,
 \cite{Weg}, imply that there exist non-trivial pairs of smooth strictly convex curves $\Gamma$
 and $\gamma$ with the desired equal tangent property, and
 moreover, the length of the tangent segments is constant for all points of $\Gamma$.
 In this setting, one can prove that if $\Gamma$ is the boundary of a body which floats in equilibrium in all positions
 and $\gamma$ (the boundary of its floating body ) is homothetic to $\Gamma$ then the curves are concentric circles.
 We do not dwell on the proof here.
}
\end{remark}

\section{In Euclidean space} \label{space}
In this section we shall see that the situation in Euclidean 3-space is very different from the plane.

Let $M$ be a $C^2$-smooth strictly convex closed surface in
$\mathbb{R}^3$. Denote by $H$ the set of points $x$ in the
exterior of $M$ such that all the tangent segments from $x$ to $M$
have equal lengths. The following theorem states that $M$ is a
sphere, provided $H$ is large enough.

\begin{theorem} \label{sphere}
Suppose that $H$ is\\
(i) a closed surface containing $M$ in its interior;\\
(ii) a plane;\\
(iii) the union of three distinct lines.\\
Then $M$ is the sphere.
\end{theorem}

\emph{Proof.} Let $x$ be a point outside of $M$. Denote by
$\gamma_x$ the curve on $M$ consisting of the contact points
between the tangents to $M$ from $x$ and $M$. Since all the
tangent segments from $x$ to $M$ have the same length,  the curve
$\gamma_{x}$ belongs to a sphere $S(x)$  centered at $x$. Hence
$\gamma_{x}$ is a line of curvature of $S(x)$. The surfaces $M$
and $S(x)$ are orthogonal along the curve $\gamma_{x}$. Therefore,
by Joachimstahl's theorem\footnote{Let two surfaces intersect
along a curve $\gamma$, and the  angle between the surfaces along
$\gamma$ is constant.  If $\gamma$ is a line of curvature on one
surface then it is also a line of curvature on the other one.},
$\gamma_{x}$ is also a line of curvature of $M$.

The idea of the proof is to show that almost every (and then, by
continuity, every) point of $M$ is umbilic. Through a non-umbilic
point there pass exactly two lines of curvature, so if one has
three such lines through a point then this point is umbilic.

To prove (i) and (ii), pick a point $p\in M$. Consider the
intersection curve of the tangent plane $T_p M$ with the surface
$H$. Choose three points $x_1,x_2,x_3$ on this curve. Then the
curves $\gamma_{x_i},\ i=1,2,3$, are different lines of curvature
on $M$ through point $p$. Hence $p$ is umbilic.

Likewise, in case (iii), let $p\in M$ be such a point that $T_p M$
intersects each of the three lines that constitute $H$ at a single
point. Almost every point $p$ satisfies this condition. Denoting
the intersection points by $x_1,x_2,x_3$, we repeat the argument
from the previous paragraph. \fin

\section{Further results} \label{further}
 The following optical (or billiard) property of ellipses is well
 known (see, e.g., \cite{Tab2}). Let $\mathcal{E}$ be an ellipse with
 the foci $P$ and $Q$, and let $X$ be a point outside of $\mathcal{E}$.
 Let $\ell_1$ and $\ell_2$ be the tangent lines to $\mathcal{E}$ from $X$.
 Then the angles between the pairs of lines $\ell_1$ and $XP$, and $\ell_2$ and $XQ$, are equal.

One has the following converse characterization of ellipses, somewhat in the spirit of Lemma \ref{line}.

\begin{proposition} \label{ellipse}
Let $\ell$ be a line tangent to a convex body $K$ in the plane, and let $P$
and $Q$ be two points in the interior of $K$. For every point $X$ in $\ell$
consider the other tangent line, $L_X$, to $K$. Suppose that the angle
between $\ell$ and $XP$ is equal to the angle between $L_X$ and $XQ$.
Then, $K$ is an ellipse with foci $P$ and $Q$. See Figure 4.
\end{proposition}

\centerline{\psset{unit=.9cm}
\begin{pspicture}(-1,-.3)(8,4)
\psellipticarc[linewidth=0.035,fillcolor=green,fillstyle=solid](3,1.5)(3,1.6){0}{180}
\pscurve[linewidth=0.035,fillcolor=green,fillstyle=solid](0,1.5)(.3,1)(5.3,.2)(6,1.5)
\psline[linecolor=blue,linewidth=.03](-1,3.1)(7.5,3.1)
\psline[linecolor=gray,linewidth=.03,linestyle=dashed](1.4,1.5)(3.7,3.1)(4.6,1.5)
\psline[linecolor=red,linewidth=.03](0,4)(7,2.3)
\psline(4.6,1.5)(4.95,2.8) \psline(1.4,1.5)(1.94,3.53)
\psline(4.6,1.5)(4.6,3.1) \psline(1.4,1.5)(1.4,3.1)
\psarc(3.7,3.1){.4}{180}{216} \psarc(3.7,3.1){.4}{300}{348}
\psdots[dotsize=3pt](1.4,1.5)(4.6,1.5)(3.7,3.1)(4.95,2.8)(1.94,3.53)(3,3.1)(1.4,3.1)(4.6,3.1)
\rput(2.4,2.8){\small $\alpha$}\rput(3.6,2.8){\small $\alpha$}
\rput(3.7,3.4){\small $X$}\rput(6,3.3){\small
$\ell$}\rput(1.2,1.3){\small $P$}\rput(4.9,1.5){\small $Q$}
\rput(6.9,2.6){\small $L_X$}\rput(.7,.1){\small $\gamma$}
\rput(1.95,3.8){\small $B$}\rput(5.2,2.94){\small $A$}
\rput(1.3,3.3){\small $R$}\rput(4.6,3.3){\small $S$}
\end{pspicture}}
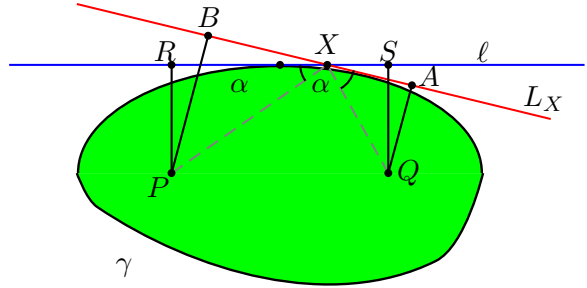
\centerline{Figure 4.}

\emph{Proof.} Suppose that the angle between $L_X$ and $XQ$ is
smaller than the angle between $L_X$ and $XP$. Let $A$ and $B$ be
the projections of $Q$ and $P$ on the line $L_X$ and let $R$ and
$S$ be the projections of $P$ and $Q$ on the line $\ell$. We
conclude from the hypothesis of the proposition that
 the right triangles $\triangle XQA$ and $\triangle XPR$
are similar, and so are the right triangles $\triangle XQS$ and
$\triangle XPB$. From these similarities we conclude that
$$RP\cdot SQ=PB\cdot QA. $$
Let $\mathcal{E}$ be the ellipse with foci $P$ and $Q$, tangent to $\ell$.
Let $L_X'$ be the tangent line to $\mathcal E$, parallel to $L_X$, such
that the ray $QA$ intersects $L_X'$. Let $A'$ and $B'$ be the
projections of $Q$ and $P$ on $L_X'$. Using the optical property of ellipses, one concludes that $RP\cdot
QS=PB'\cdot QA'$. It follows that  $L_X'$ coincides with $L_X$. Thus
the tangent lines to $\mathcal{E}$ and $K$ from all points of the line $\ell$ coincide, and hence $K=\mathcal{E}$.
 \fin

In conclusion, we remark that Lemma \ref{line} holds in all three classical geometries: elliptic,
 Euclidean, and hyperbolic;  we give a proof that works in all three cases, cf. \cite{Tab1,Tab3}.

In the argument below, a ``circle" means a curve of constant curvature. In Euclidean and
elliptic geometry this is indeed a circle; in hyperbolic geometry this may be a circle, a
horocycle, or an arc with both endpoints at infinity, depending on the value of the curvature.

Let $x$ be a point of $\ell=T_p \gamma$, and let $xy,\ y\in\gamma$, be the other tangent segment
to $\gamma$ from $x$. Then $|xy|=|xp|$ if and only if there exists a ``circle"  tangent to $\gamma$ at $p$ and $y$.

Consider the family of ``circles", tangent to $\ell$ at point $p$. In the complement of $p$,
these curves form a smooth foliation $\cal F$. Since the two tangent segments to $\gamma$ from
every point $x\in\ell$ have equal lengths, the curve $\gamma$ is everywhere tangent to the leaves
of the foliation $\cal F$. It follows that $\gamma$ coincides with a leaf, that is, $\gamma$ is
a ``circle". Since $\gamma$ is a closed curve, it is indeed a circle.

\medskip

{\bf Acknowledgments}.  S. T. was partially supported by the Simons Foundation grant No 209361 and by the NSF grant DMS-1105442.

\end{document}